\magnification=1200%
\emergencystretch = 10 pt%
\hfuzz            = 10 pt%
\font\bb=msbm10%

\def\square{\hbox{\rlap{$\sqcap$}$\sqcup$}}

\centerline{\bf ON PTOLEMAIC METRIC SIMPLICIAL COMPLEXES}
\bigskip\bigskip
\centerline{\bf S. M. Buckley, J. McDougall, D. J. Wraith}
\bigskip\bigskip
\itemitem{} {\bf Abstract:} We show that under certain mild conditions, a 
metric simplicial
complex which satisfies the Ptolemy inequality is a CAT(0) space. 
Ptolemy's inequality is closely related to inversions of metric
spaces. For a large class of metric simplicial complexes, we
characterize those which are isometric to Euclidean space in terms 
of metric inversions.
\bigskip\bigskip

\centerline{\bf \S1 Introduction}
\bigskip

The aim of this paper is to study some geometric aspects of metric
simplicial complexes. (See \S2 for definitions, or the book [BH] for a
comprehensive introduction.) Roughly speaking, in a metric simplicial
complex the simplexes are all subsets of a fixed Riemannian manifold with
constant sectional curvature, and are glued together by isometries. There
is a natural way to define a (distance) metric on such an object, and
consequently these objects form a large and interesting class of metric
spaces.

We are specifically interested in those metric simplicial complexes $(K,d)$
which are {\it Ptolemaic spaces}.
This means that the following {\it Ptolemy inequality} is satisfied for
every quadruple of points $x,y,z,p \in (K,d):$
$$ d(x,y)d(z,p) \le d(x,z)d(p,y) + d(x,p)d(y,z). $$
The significance of this inequality in various metric space settings has
been studied recently, for example in [BFW] and [FLS]. Note that a
classical result says that Ptolemy's inequality holds in the Euclidean
plane, with equality if and only if the points $x$, $y$, $z$, and $p$ lie
on a circle in that order.

The Ptolemy inequality is closely related to the concept of {\it metric
space inversion}. Inversion (or reflection) about the Euclidean unit sphere
is a bijection on $\hbox{\bb R}^n \backslash \{0\}$, so we can pull back
Euclidean distance to get a new distance on $\hbox{\bb R}^n \backslash
\{0\},$ namely $i_0(x,y):=|x-y|/|x|\,|y|$. Inversion has been generalized
to the setting of a metric space $(X,d)$ in [BHX]: for fixed $p\in X$,
define
$$
i_p(x,y)={{d(x,y)} \over {d(x,p)d(y,p)}}\,, \qquad x,y \in X_p\,,
$$
where $X_p:=X\backslash \{p\}$. In general this is not a metric on $X_p$,
but a metric $d_p:X_p\times X_p\to[0,\infty)$ can be defined which is both
subordinate to, and comparable with, $i_p$; see [BHX; Lemma 3.2].

The definition of $d_p$ is more complicated than that of $i_p$, so it is
natural to ask when $i_p$ itself is a metric for all $p\in X$. This reduces
to deciding if $i_p$ satisfies the triangle inequality, and it is
elementary to observe that $(X_p,i_p)$ is a metric space for all $p\in X$
if and only if $(X,d)$ is Ptolemaic.

It is not hard to see that CAT(0) spaces are Ptolemaic; see for example
[BFW; \S3]. As a corollary, the inversions $i_p$ are metrics in every
CAT(0) space. On the other hand, the converse statement is not true: there
exists a geodesic Ptolemaic space that is not CAT(0); see the comments
after Theorem 1.1 in [FLS].

Nevertheless, if more structure is imposed on the space under
consideration, it is possible to give converse statements. In particular, a
complete Riemannian or Finsler manifold $M$ is Ptolemaic if and only if it
is CAT(0), or equivalently a Hadamard manifold; for a proof, see [BFW] or
[K].

The first result in this paper shows that an analogous result is true for
metric simplicial complexes.

\proclaim Theorem A. Let $K$ be a metric simplicial complex with simplexes
of curvature $\kappa$, dimension $n \ge 2$ and $\hbox{Shapes}(K)$ finite.
If $K$ satisfies the Ptolemy inequality, then $K$ must be CAT(0). In
particular we must have $\kappa \le 0.$ \par

The condition $\lq\hbox{Shapes}(K)$ finite' above means that the complex
contains only finitely many isometry types of simplex.

In [BFW, Section 6], it was also shown that inversion can be used to
characterize Euclidean space amongst all Riemannian manifolds. In fact if
$M$ is a complete Riemannian manifold, then the inversion of $M$ with
respect to $p$ has the structure of a Riemannian manifold for all $p\in M$
if and only if $M$ is Euclidean space. Moreover, the inversion is a
manifold if and only if it is a length space. Our second main result is an
analogue of this result for metric simplicial complexes.

\proclaim Theorem B. Let $K$ be a metric simplicial complex with
$\hbox{Shapes}(K)$ finite, which is homeomorphic to $\hbox{\bb R}^n.$ If
for all $p\in K$, $i_p$ is a metric and the inversion of $K$ with respect
to $p$ is a length space, then $K$ must be isometric to $\hbox{\bb R}^n.$
\par

Note that some topological assumption such as the condition that $K$ is
homeomorphic to $\hbox{\bb R}^n$ is required in Theorem B, since otherwise
there are some trivial counterexamples, such as a complex consisting of a
single simplex, or certain complexes containing simplexes of differing
dimensions.

The remaining sections of this paper are laid out as follows. In \S2 we
give all the definitions and background results that we will need. In \S3
we will prove Theorem A and in \S4 we will prove Theorem B.
\bigskip

\centerline{\bf \S2 Definitions and background results}
\bigskip

We begin this section by recalling some special types of metric space.

We say that $(X,d)$ is a {\it length space} if the distance $d(x,y)$
between any pair of points is always equal to the infimum of the lengths of
paths between the points. (See [BBI; Chapter 2] a full description of the
notion of path length in a metric space setting.) A path $\gamma$ of length
$d(x,y)$ joining $x,y\in X$ is called a {\it geodesic segment}, and is
often denoted $[x,y]$. We call $(X,d)$ a {\it geodesic space} if all pairs
of points can be joined by geodesic segments, that is, the above infimum is
always attained.

A {\it geodesic triangle} $T(x,y,z)$ is a collection of three points
$x,y,z\in X$ together with a choice of geodesic segments $[x,y]$, $[x,z]$
and $[y,z]$. Given such a geodesic triangle $T(x,y,z)$, a {\it comparison
triangle} will mean a geodesic triangle in a simply-connected constant
curvature surface $\bar{T}(\bar{x},\bar{y},\bar{z})$, such that
corresponding distances coincide: $d(x,y)=\bar{d}(\bar{x},\bar{y})$,
$d(y,z)=\bar{d}(\bar{y},\bar{z})$, $d(z,x)=\bar{d}(\bar{z},\bar{x})$.
(Usually, such a comparison triangle will be in the Euclidean plane.)

Recall the definition of a $\hbox{CAT}(\kappa)$ space. This is a geodesic
metric space $(X,d)$ with the following property. Let $T$ be a geodesic
triangle in $X$, and $\bar{T}$ a comparison triangle in constant curvature
$\kappa.$ Let $D_\kappa$ denote the diameter of the unique simply-connected
surface of constant curvature $\kappa$ (so $D_\kappa=\infty$ if $\kappa \le
0$). If the perimeter of $T$ is less than $2D_\kappa$, then given any two
points $x,y \in T$ and corresponding comparison points $\bar{x}, \bar{y}
\in \bar{T}$ we have $$d(x,y) \le \bar{d}(\bar{x},\bar{y}).$$ We say that
$T$ is {\it equal to} $\bar{T}$ if this inequality is actually an equality
for all pairs of points $x,y.$

A metric space is said to be {\it proper} if all its closed balls are
compact. \proclaim Lemma 2.1. A length space is proper if and only if it is
locally compact and complete. \par \noindent {\bf Proof.} That locally
compact and complete imply proper for length spaces is precisely [BBI;
2.5.22]. In the other direction the local compactness is trivial and the
completeness follows from the fact that any Cauchy sequence can be
contained in a closed ball. \hfill{\square}
\smallskip

The {\it one-point extension of $X$} is defined to be
$$
\hat{X}:=\left \{\matrix{ X \hbox{ when } X \hbox{ is bounded}\,, \cr\cr
X\cup\{\infty\} \hbox{ when } X \hbox{ is unbounded}\,;
} \right.
$$
the open sets in $\hat{X}$ include those in $X$ together with complements
(in $\hat{X}$) of closed balls (in $X$).  Thus when $X$ is a proper space,
$\hat{X}$ is simply its one-point compactification.

Suppose now that $i_p$ is actually a metric on $X$. Recall that this is the
case if $(X,d)$ is Ptolemaic. Additionally, it is worth noting that $i_p$
is a metric if and only if $i_p=d_p$, where $d_p$ is the metric defined in
[BHX] mentioned in the introduction.

When $(X,d)$ is unbounded, there is a unique point $p'$ in the completion
$(\hat{X}_p,\hat{i}_p)$ of $(X_p,i_p)$ which corresponds to the point
$\infty$ in $\hat{X}$. (Any unbounded sequence in $(X_p,d)$ is a Cauchy
sequence in $(X_p,i_p)$, and any two such sequences are equivalent.) Note
that $\hat{i}_p(x,p')=1/d(x,p).$ We denote this completion
$\hbox{Inv}_p(X)$ and refer to it as the {\it inversion\/} of $(X,d)$ with
respect to the base point $p$. For example, with this definition,
$\hbox{Inv}_p(X)$ will be complete (or proper) whenever $X$ is complete (or
proper).
\smallskip

We now turn our attention to metric simplicial complexes. We follow [BH;
p.~98]. Let $M^n_{\kappa}$ denote the simply connected $n$-manifold with
constant (sectional) curvature $\kappa.$ If $n \le m$ then an $n$-plane in
$M^m_{\kappa}$ will be a subspace isometric to $M^n_{\kappa}.$ If $(n+1)$
points in $M^m_{\kappa}$ do not lie in any $(n-1)$-plane, we say the points
are in {\it general position}. A {\it geodesic $n$-simplex} in
$M^m_{\kappa}$ is defined to be the convex hull of $(n+1)$ points in
general position. If $\kappa>0$ then the vertices (that is, the points
defining the convex hull) must lie in an open ball of radius
$D_{\kappa}/2,$ with $D_{\kappa}$ as above.

\proclaim Definition 2.2. Let $\{S_{\lambda} \,|\, \lambda \in \Lambda\}$
be a collection of geodesic simplexes with $S_{\lambda} \subset
M^{n_{\lambda}}_{\kappa}.$ Let $X=\cup_{\lambda \in \Lambda}(S_{\lambda}
\times \{\lambda\}),$ let $\sim$ be an equivalence relation on $X$ and set
$K:=X/\sim.$ Let $p:X \rightarrow K$ be the quotient map and set
$p_{\lambda}(x):=p(x,\lambda).$ $K$ is called an {\it
$M_{\kappa}$-simplicial complex} if the following two conditions are
satisfied:
\item{(1)} $p_{\lambda}$ is injective for all $\lambda \in \Lambda$;
\smallskip
\item{(2)} if $p_{\lambda}(S_{\lambda}) \cap p_{\lambda'}(S_{\lambda'})
    \neq \emptyset$ then there is an isometry $h_{\lambda,\lambda'}$ from a
    face $T_{\lambda} \subset S_{\lambda}$ to a face $T_{\lambda'} \subset
    S_{\lambda'}$ such that $p(x,\lambda)=p(x',\lambda')$ if and only if
    $x'=h_{\lambda,\lambda'}(x).$
\par
In this paper, we will always assume that our complexes are connected.
\proclaim Definition 2.3. The set of isometry classes of the faces of the
geodesic simplexes $S_{\lambda}$ will be denoted $\hbox{Shapes}(K).$ \par

An $M_{\kappa}$-complex $K$ becomes a pseudometric space when equipped with
the following natural metric. Given any two points $x$ and $y$ in $K,$
consider a sequence of points $x=x_1,x_2,...,x_n=y$ with the property that
every adjacent pair of points belong to a common simplex. As a result, the
distance between any two adjacent points can be found, and so the `length'
of the sequence can be computed. The pseudometric $d(x,y)$ is then defined
to be the infimum of the lengths of all such sequences linking $x$ and $y.$
If this pseudometric is actually a metric, we call $K$ a {\it metric
simplicial complex.}

The following result from Bridson's thesis (see [BH; p.~97]) shows that
metric simplicial complexes have good properties:

\proclaim Theorem 2.4.
An $M_{\kappa}$-simplicial complex $K$ with $\hbox{Shapes}(K)$ finite is a
metric simplicial complex, and moreover it is a complete geodesic space.


With a view towards exploring some of the more detailed structure of metric
simplicial complexes, we introduce the concept of a $\kappa$-cone over a
metric space [BH; p.~59]. Given $\kappa \in \hbox{\bb R}$ and a metric
space $(X,d)$, the $\kappa$-cone over $X$, $C_{\kappa}X$ is given by $X
\times [0,\infty)/\sim$ if $\kappa \le 0$ and $X \times
[0,D_{\kappa}/2]/\sim$ in the case $\kappa>0,$ where in either case $(t,x)
\sim (t',x')$ if and only if $t=t'=0$ or we have equality of pairs. The
equivalence class corresponding to the point 0 is the {\it vertex} of the
cone. We define a metric $d_C$ on $C_{\kappa}X$ as follows. Let $p=[x,t]$
and $p'=[x',t'].$
\item{} If $\kappa=0$ then set
$$
d_C(p,p')=t^2+t'^2-2tt'\cos(\min\{\pi,d(x,x')\}).
$$
\item{} If $\kappa<0$ then define $d_C$ via
$$
\eqalign{
  \cosh(\sqrt{-\kappa}) &d_C(p,p') = \cr
 &\cosh(\sqrt{-\kappa}t)\cosh(\sqrt{-\kappa}t') -
  \sinh(\sqrt{-\kappa}t)\sinh(\sqrt{-\kappa}t')\cos(\min\{\pi,d(x,x')\}). }
$$
\item{} If $\kappa>0$ then define $d_C$ via
$$
\cos(\sqrt{\kappa}) d_C(p,p') =
  \cos(\sqrt{\kappa}t)\cos(\sqrt{\kappa}t') +
  \sin(\sqrt{\kappa}t)\sin(\sqrt{\kappa}t')\cos(\min\{\pi,d(x,x')\}).
$$

\proclaim Definition 2.5. For any point $x$ in a geodesic simplex $S$, the
{\it link} of $x$ in $S$, $\hbox{Lk}(x,S)$ is the set of unit vectors at
$x$ which point into $S$. If $x$ in an element of a metric simplicial
complex $K$, the link of $x$ in $K$, $\hbox{Lk}(x,K)$ is the union of the
links of $x$ in all simplexes to which $x$ belongs. (See [BH; pp.~102-103]
for an alternative description.) \par

One can define a natural (angular) pseudometric on $\hbox{Lk}(x,K)$.

\proclaim Fact 2.6. {$\!$\rm [BH; p.~103] }
For any $M_{\kappa}$-complex $K$ and any point $x \in K$, $\hbox{Lk}(x,K)$
is an $M_1$-simplicial complex. Hence if $\hbox{Shapes}(K)$ is finite, then
$\hbox{Lk}(x,K)$ is a metric simplicial complex and a complete geodesic
space.

Note that an $M_{\kappa}$-complex is not necessarily a $\hbox{CAT}(\kappa)$
space. In particular, a link complex is not necessarily a $\hbox{CAT}(1)$
space. This motivates the following definition: \proclaim Definition 2.7.
An $M_{\kappa}$-simplicial complex $K$ satisfies the {\it link condition}
if for every vertex $v \in K,$ the link complex $\hbox{Lk}(v,K)$ is a
$\hbox{CAT}(1)$ space. \par

The following result will be crucial in the proof of Theorem A in \S2. It
is a combination of [BH; II.5.4] and [BH; II5.6].

\proclaim Theorem 2.8.
Let $K$ be an $M_{\kappa}$-simplicial complex with $\hbox{Shapes}(K)$
finite.
\item{(a)} If $\kappa \le 0$, then the following conditions are equivalent:
\smallskip
\itemitem{i)} $K$ is a $\hbox{CAT}(\kappa)$ space;
\smallskip
\itemitem{ii)} $K$ satisfies the link condition and contains no
    isometrically embedded circles.
\smallskip
\item{(b)} If $\kappa > 0$, then the following conditions are equivalent:
\smallskip
\itemitem{i)} $K$ is a $\hbox{CAT}(\kappa)$ space;
\smallskip
\itemitem{ii)} $K$ satisfies the link condition and contains no
    isometrically embedded circles of length less than
    $2\pi/\sqrt{\kappa}$.
\item{(c)} If $K$ is a two-dimensional complex, then it satisfies the link
    condition if and only if for each vertex $v \in K,$ every injective
    loop in $\hbox{Lk}(v,K)$ has length at least $2\pi.$
\par

The final result in this section is essentially [BH; I.7.17] in the special
case of a complex $K$ with $\hbox{Shapes}(K)$ finite.

\proclaim Theorem 2.9.
For an $M_{\kappa}$-simplicial complex $K$ and a point $x \in K$, the
$\epsilon$-ball about $x$, $B(x,\epsilon),$ is isometric to the
$\epsilon$-ball about the vertex in $C_{\kappa}(\hbox{Lk}(x,K))$ for all
$\epsilon$ sufficiently small. \par
\bigskip\bigskip

\centerline{\bf \S3 The Proof of Theorem A}
\bigskip\bigskip

Recall Theorem A from the Introduction:%

\proclaim Theorem A.
Let $K$ be an $M_{\kappa}$-simplicial complex of
dimension $n \ge 2$ with $\hbox{Shapes}(K)$ finite.  If $K$ satisfies the
Ptolemy inequality, then $K$ must be CAT(0). \par
\medskip

\noindent {\bf Proof.}
First note that the Ptolemaic condition must hold in each simplex of $K$,
so as noted in the Introduction this means each simplex must individually
be a $\hbox{CAT}(0)$ space. In turn this means that we must have $\kappa
\le 0$.

As $K$ is Ptolemaic, it cannot contain an isometrically embedded circle, as
the $\lq$quarter points' would violate the Ptolemy inequality. Therefore,
by Theorem 2.8, we see that $K$ is CAT($\kappa$) if and only if it
satisfies the link condition at each vertex. (The $\hbox{Shapes}(K)$
condition is required for this theorem.) By the same result, the link
$Lk(v,K)$ of some vertex $v$ of $K$ is a CAT(1) space if and only if
$Lk(v,K)$ satisfies the link condition at each of {\it its} vertices and
contains no isometrically embedded $S^1$ of length strictly less than
$2\pi$. We will use the term $\lq$short' to describe a circle of length
less than $2\pi$. Therefore $K$ fails to satisfy the link condition if and
only if the link in $Lk(v,K)$ of some vertex either fails to satisfy the
link condition or $Lk(v,K)$ has a short isometrically embedded circle.

If the link condition fails in $Lk(v,K)$ then either the link condition
fails in a link of the link, or some link of the link has a short
isometrically embedded $S^1$. (From now on we will suppress the vertices
from the link notation, so $Lk(K)$ will denote the original link, a link of
this link will be written $Lk^2(K)$ and so on.) If the link condition keeps
failing in $Lk(K)$, $Lk^2(K)$, $Lk^3(K)$,... then eventually the link
condition will fail in a 1-complex. But the link condition fails in a
1-complex if and only if there is an isometrically embedded (i.e.
injective) short loop. We therefore arrive at the following

\proclaim Observation.
A Ptolemaic $M_{\kappa}$-complex $K$ fails to be a CAT($\kappa$) space
($\kappa \le 0$) if and only if some $Lk^m(K)$, $m \ge 1$, has an
isometrically embedded short loop. \par

Our strategy is to show that the existence of an isometrically embedded
$S^1$ in some $Lk^m(K)$ gives a contradiction with the original complex $K$
being Ptolemy. More specifically, the quarter-points in this $S^1$ can be
associated to four points in $K$ (by the {\it $\lq$method of association'}
below), for which the Ptolemy inequality can be shown to fail. Thus $K$ has
to be CAT($\kappa$), and in particular CAT(0).
\smallskip

\noindent {\it Method of Association.} To the $S^1$ in $Lk^m(K)$ we
actually associate an $S^1$ in $K$. The desired four points are then the
points in $K$ corresponding to the quarter-points in the original circle.
\par

There is an $\epsilon>0$ such that the $\epsilon$-neighbourhood of the cone
point in $C_1 Lk^m(K)$ is isometric to the $\epsilon$-neighbourhood of the
central vertex for $Lk^{m-1}(K)$. (Here, $Lk^0(K)$ should be interpreted as
$K$ itself.) The $S^1 \subset Lk^m(K)$ gives a topological circle $\epsilon
S^1 \subset C_1 Lk^m(K)$, and in turn this gives a topological circle in
$Lk^{m-1}(K)$. Continuing in this way (with the same suitably small
$\epsilon$) we will eventually produce an embedded $S^1$ in $K$.

In order to show that the Ptolemy inequality fails for our four points in
$K$, we need to investigate how distance alters when we embed points into
successive 1-cones, and ultimately into $\kappa$-cones, with $\kappa \le
0$.

Consider two points separated by a distance $D$ in some metric space $Y$.
(Assume $D \le \pi$.) Then according to the cone metric definitions given
in \S2, the separation of the corresponding points at a distance $\epsilon$
from the vertex in $C_1 Y$ is
$$
\eqalign{&\cos^{-1}(\cos^2 \epsilon +(\sin^2 \epsilon) \cos D) \cr
=&\cos^{-1}(1-(\sin^2 \epsilon)(1-\cos D)). \cr}
$$
Similarly, the separation of the corresponding points at a distance
$\epsilon$ from the vertex in the $\lq$cone over the cone' $C_1^2 Y$ is
$$
\eqalign{
&\cos^{-1}(1-(\sin^2 \epsilon)
  [1-\cos(\cos^{-1}[(1-\sin^2 \epsilon)(1-\cos D)])]) \cr =
&\cos^{-1}(1-(\sin^2 \epsilon)[1-(1-(\sin^2 \epsilon)(1-\cos D))]) \cr
=& \cos^{-1}(1-(\sin^4 \epsilon) (1-\cos D)).\cr}
$$
Similarly, the corresponding distance in $C_1^r Y$ is
$\cos^{-1}(1-(\sin^{2r} \epsilon) (1-\cos D)).$
\medskip

\noindent {\bf The case $\kappa=0$.} \smallskip

Suppose the points at the above separation are finally embedded at a
distance $\epsilon$ from the vertex of a 0-cone. By the cone metric
definitions in \S2, the separation is then
$$
\eqalignno{
[2\epsilon^2 &- 2\epsilon^2
  \cos(\cos^{-1}[1-(\sin^{2r} \epsilon)(1-\cos D)])]^{1 \over 2}\cr
&=\, \epsilon\sqrt 2 [1-(1-(\sin^{2r} \epsilon)(1-\cos D))]^{1 \over 2} \cr
&=\, \epsilon \sqrt 2 (\sin^{r} \epsilon) \sqrt{1-\cos D}. &(\ast) }
$$
In the original isometrically embedded $S^1$, suppose the separation of
adjacent quarter points is $D<\pi/2$. The separation of opposite points is
then $2D$.
\par
We check Ptolemy's inequality for the corresponding points in $K$. Using
($\ast$) we see that Ptolemy will fail if
$$
2\epsilon^2 (\sin^{2r} \epsilon) (1-\cos 2D)>
  4\epsilon^2 (\sin^{2r} \epsilon) (1-\cos D);
$$
that is, if
$$ 1-\cos 2D >2(1-\cos D). $$
But $\cos 2D=2\cos^2 D -1$, so this inequality is really
$$ 2-2\cos^2 D>2-2\cos D, $$
or more simply
$$\cos D>\cos^2D.$$
But $D \in(0,\pi/2)$, therefore $\cos D \in (0,1)$ and so the inequality
must be true.
\medskip

\noindent {\bf The case $\kappa<0$.} \smallskip

As in the $\kappa=0$ case, suppose the four points with (adjacent)
separation $\cos^{-1}(1-(\sin^{2r} \epsilon) (1-\cos D))$ are finally
embedded at a distance $\epsilon$ from the vertex of a $\kappa$-cone, with
$\kappa<0.$ By the cone metric definitions in \S2, the separation is then
$$
\cosh^{-1}(\cosh^2 \epsilon - (\sinh^2 \epsilon) \cos
(\cos^{-1}[1-(\sin^{2r}\epsilon) (1-\cos D)]))
$$
$$
=\cosh^{-1}(\cosh^2 \epsilon - (\sinh^2 \epsilon)
  [1-(\sin^{2r}\epsilon) (1-\cos D)]),
$$
and the separation of opposite points is
$$
=\cosh^{-1}(\cosh^2 \epsilon - (\sinh^2 \epsilon)
  [1-(\sin^{2r}\epsilon) (1-\cos 2D)]).
$$
We show that for $\epsilon$ sufficiently small, Ptolemy's inequality fails
for these distances.

Setting $D'=\cos^{-1}(1 - (\sin^{2r} \epsilon) (1-\cos D))$ and
$D''=\cos^{-1}(1 - (\sin^{2r} \epsilon) (1-\cos 2D)),$ let us label the
corresponding separations $S'(\epsilon)$ and $S''(\epsilon)$ respectively.
It follows from the cone metric definitions that these separations are
equal to the length of the third side in an isosceles triangle in the
simply connected space of constant curvature $\kappa$, where the equal
sides have length $\epsilon,$ and the angle between the sides is $D'$
respectively $D''$.

It will be convenient to find alternative expressions for $S'(\epsilon)$
and $S''(\epsilon)$ based on the sine law for triangles in hyperbolic
space. Note that for a geodesic triangle in the simply connected space of
constant curvature $\kappa<0,$ the sine law reads $${{\sinh a} \over {\sin
A}}= {{\sinh b} \over {\sin B}}={{\sinh c} \over {\sin C}}$$ where $a,b,c$
are the side lengths and $A,B,C$ are the angles (see [C; p.~94]).

Consider splitting each isosceles triangle into two equal triangles by
introducing a line dividing the angle $D'$ (respectively $D''$) in half.
The point at which this line meets the opposite side is clearly the point
on that side closest to the cone vertex. The angle made between the two
lines is therefore $\pi/2.$

Applying the sine law to one of our $\lq$half-triangles' gives
$$
{1 \over {\sqrt{-\kappa}}}
  {{\sinh (\epsilon\sqrt{-\kappa})} \over {\sin (\pi/2)}}
= {1 \over {\sqrt{-\kappa}}}
  {{\sinh (\sqrt{-\kappa}S'(\epsilon)/2)} \over {\sin (D'/2)}}\, ,
$$
which after rearranging gives
$$
S'(\epsilon) = {2 \over {\sqrt{-\kappa}}}
  \sinh^{-1}(\sin (D'/2)\sinh (\epsilon\sqrt{-\kappa})).
$$
Similarly for $S''$ and $D''$. Expanding this as a Taylor series about
$\epsilon=0$ gives
$$
\eqalign{
S'(\epsilon)  &= 2\epsilon \sin (D'/2)+ O(\epsilon^3), \cr
S''(\epsilon) &= 2\epsilon \sin (D''/2)+ O(\epsilon^3). }
$$

The Ptolemy inequality will fail for our chosen points if
${S''}^2(\epsilon)>2{S'}^2(\epsilon),$ that is, if
$$
4\epsilon^2 \sin^2 (D''/2) + O(\epsilon^4) >
  8\epsilon^2 \sin^2(D'/2) + O(\epsilon^4),
$$
or equivalently
$$\sin^2(D''/2)+O(\epsilon^2)>2\sin(D'/2)+O(\epsilon^2). \eqno{(\dag)}$$
We claim that this is true for all sufficiently small $\epsilon.$ Before we
can establish this, however, we need two lemmas.

\proclaim Lemma 3.1.
For $\lambda>\sqrt 2$ and $x$ sufficiently small, the following inequality
holds:
$$ \sin^2 x > 2\sin^2(x/\lambda). $$

\noindent {\bf Proof.}
For $x$ small we have $\sin^2 x=x^2+O(x^4).$ Therefore
$$ 2\sin^2(x/\lambda)=2(x/\lambda)^2+O(x^4). $$
For $x$ so small that the $O(x^4)$ terms are irrelevant, establishing the
inequality reduces to showing that $x^2>2(x/\lambda)^2,$ which is true
since $\lambda > \sqrt 2.$ \hfill{\square}

\proclaim Lemma 3.2.
Given $D \in (0,\pi/2)$, choose $\lambda \in (\sqrt 2,\sqrt{2(1+\cos D)}).$
Then for $\epsilon>0$ sufficiently small (depending on $\lambda$) we have
$$ \cos^{-1}(1-\epsilon (1-\cos 2D))> \lambda \cos^{-1}(1-\epsilon (1-\cos
D)).$$

Notice that both sides would be zero if we were allowed to set
$\epsilon=0.$ We examine the derivatives with respect to $\epsilon$ of each
side in the above inequality.

$$
{d \over {d\epsilon}}\cos^{-1}(1-\epsilon (1-\cos 2D)) =
 {{1-\cos 2D} \over {\sqrt{2\epsilon(1-\cos 2D)-
 \epsilon^2(1-\cos 2D)^2}}}\,;
$$
$$
{d \over {d\epsilon}}\cos^{-1}(1-\epsilon (1-\cos D)) =
 {{1-\cos D} \over {\sqrt{2\epsilon(1-\cos D)-
 \epsilon^2(1-\cos D)^2}}}\,.
$$
We can therefore establish the truth of our inequality by showing that for $\epsilon$ sufficiently small:
$$
{{1-\cos 2D} \over {\sqrt{2\epsilon(1-\cos 2D)-\epsilon^2(1-\cos 2D)^2}}} >
 {{\lambda(1-\cos D)} \over {\sqrt{2\epsilon(1-\cos D)-
 \epsilon^2(1-\cos D)^2}}}\,;
$$
or equivalently
$$
{{1-\cos 2D} \over {\lambda(1-\cos D)}} >
 \sqrt{{{2(1-\cos 2D)-\epsilon(1-\cos 2D)^2} \over {2(1-\cos D)-\epsilon
 (1-\cos D)^2}}}.
$$
As $\epsilon \rightarrow 0$, the right hand side tends to
$$ \sqrt{{{1-\cos 2D} \over {1-\cos D}}}. $$
It therefore suffices to show that
$$
{{1-\cos 2D} \over {\lambda(1-\cos D)}} >
\sqrt{{{1-\cos 2D} \over {1-\cos D}}}\,;
$$
or equivalently
$$ {{1-\cos 2D} \over {1-\cos D}}>\lambda^2. $$
Using the double angle formula for $\cos 2D$ in the above and rearranging,
we obtain
$$-2\cos^2 D+\lambda^2\cos D +(2-\lambda^2)>0.$$
The roots of this quadratic expression are $\cos D=1$ and $\cos
D=(\lambda^2-2)/2.$ Therefore the quadratic expression is positive
precisely when $\cos D \in ( (\lambda^2-2)/2,1),$ assuming $\lambda<2.$ But
$\cos D<1$ anyway since $D\in(0,\pi/2)$, so we only require $\cos D
> (\lambda^2-2)/2$, i.e.~$\lambda<\sqrt{2(1+\cos D)},$ as claimed.
\hfill{\square}
\medskip

\noindent {\it Proof of Theorem A continued.}
\smallskip

Now let us return to the inequality $(\dag)$. Lemma 3.2 shows that for a
suitable choice of $\lambda$, $D''>\lambda D'$ for $\epsilon$ suitably
small. By Lemma 3.1, for this $\lambda$ and $\epsilon$ sufficiently small
we have
$$
\sin^2(D''/2) > \sin^2(\lambda D'/2)
  > 2\sin^2\Bigl({{\lambda D'} \over {2\lambda}}\Bigr) = 2\sin^2(D'/2).
$$
Thus for $\epsilon$ sufficiently small we see that $(\dag)$ holds, and thus
the Ptolemy inequality fails.
\par

Finally, note that in both the $\kappa=0$ and $\kappa<0$ cases, we did not
consider the situation where the isometrically embedded $S^1$ is in $K$, as
opposed to some $Lk^m(K).$ However this situation is trivial, as the
failure of the Ptolemy inequality is equivalent to showing $4D^2>2D^2$,
which is clearly true. \hfill{\square}
\bigskip\bigskip\bigskip

\centerline{\bf \S4 The Proof of Theorem B}
\bigskip\bigskip

First, let us recall Theorem B from the Introduction:

\proclaim Theorem B. Let $K$ be a metric simplicial complex with
$\hbox{Shapes}(K)$ finite, which is homeomorphic to $\hbox{\bb R}^n,\, n\ge
2.$ If, for all $p\in K$, $i_p$ is a metric and the inversion of $K$ with
respect to $p$ is a length space, then $K$ must be isometric to $\hbox{\bb
R}^n.$ \par

Before proving this we need a sequence of lemmas.

\proclaim Lemma 4.1.
If $(X_p,i_p)$ is a length space, then so is $\hbox{Inv}_p(X) =
(\hat{X}_p,\hat{i}_p).$
\medskip

\noindent {\bf Proof.}
It is clear that the only issue is with distances to the point $p' \in
\hbox{Inv}_p(X)$ corresponding to the point at infinity in the completion
of $(X,d)$ (in the case that $(X,d)$ is unbounded.

Consider any $d$-unbounded Cauchy sequence $\{y_i\} \subset
\hbox{Inv}_p(X).$  Given any $x \in X$ and $\epsilon>0,$ we construct a
path from $x$ to $p'$ with length strictly less than
$\hat{i}_p(x,p')+\epsilon.$ We do this as follows. Choose $I_0$ such that
for all $i,j \ge I_0$ we have $\hat{i}_p(y_i,y_j)<\epsilon/8.$ By removing
points of the sequence and re-labelling if necessary, we can assume without
loss of generality that for all $i \ge I_0,$
$$\hat{i}_p(y_i,y_{i+1})<{{\epsilon} \over 8}.2^{-i}.$$
Moreover, since $(X_p,i_p)$ is a length space, we can choose a path from
$y_i$ to $y_{i+1}$ for each $i \ge I_0$ with length strictly less than
$$\hat{i}_p(y_i,y_{i+1})+{{\epsilon} \over 8}.2^{-i}<{{\epsilon} \over 4}.2^{-i}.$$
We can also choose a path from $x$ to $y_{I_0}$ of length strictly less
than $$\hat{i}_p(x,y_{I_0})+{{\epsilon} \over 4}.$$ Therefore concatenating
this with the paths between the $\{y_i\}$, $i \ge I_0,$ gives a path of
total length at most $\hat{i}_p(x,y_{I_0})+(\epsilon/2).$ Parameterising
this path by arclength and calling it $\gamma(t)$ gives a map
$$\gamma:[0,L) \rightarrow \hat{X}_p,$$ and clearly $$\lim_{t \rightarrow
L} \hat{i}_p(\gamma(t),p')=0.$$ We can therefore $\lq$complete' the path by
adding the point $p',$ to get a path $\hat{\gamma}(t)$ defined on the
interval $[0,L].$ Obviously, adding this point does not affect the length.
We have therefore constructed a path $\hat{\gamma}$ from $x$ to $p'$ with
length at most $\hat{i}_p(x,y_{I_0})+(\epsilon/2).$ Using the triangle
inequality we see that
$$
\hat{i}_p(x,y_{I_0}) \le \hat{i}_p(x,p')+\hat{i}_p(y_{I_0},p') \le
  \hat{i}_p(x,p')+ {{\epsilon}\over 8}.
$$
Therefore the length of $\hat{\gamma}$ satisfies
$$
\hbox{length of }\hat{\gamma}
  < \hat{i}_p(x,p')+{{\epsilon} \over 2}+{{\epsilon} \over 8}
  < \hat{i}_p(x,p')+\epsilon
$$
as claimed. \hfill{\square}
\medskip

\proclaim Lemma 4.2.
Let $K$ be an $M_{\kappa}$-complex with $\hbox{Shapes}(K)$ finite, which is
homeomorphic to $\hbox{\bb R}^n,$ for which $i_p$ is a length metric. Then
$\hat{i}_p$ is a geodesic metric. \par
\medskip

\noindent {\bf Proof.}
By Lemma 4.1, $\hat{i}_p$ is also a length metric. Now a complete locally
compact length space is a geodesic space (see [BBI; 2.5.23]), and so it
suffices to show that $\hbox{Inv}_p(K)$ is complete and locally compact. By
Lemma 2.1, this in turn is equivalent to properness. However, as noted in
\S2 (or see [BHX; p.~6]), if $(X,d)$ is proper, so is $\hbox{Inv}_p(X).$
Therefore in our case, we will be done if we can show that the complex
$(K,d)$ is proper, or equivalently that it is complete and locally compact.
As local compactness is preserved by homeomorphism, and $K \cong \hbox{\bb
R}^n,$ the local compactness of $K$ is clear. Completeness is given by
Theorem 2.4. \hfill{\square}
\medskip

\proclaim Lemma 4.3.
Suppose the metric simplicial complex $K$ is such that $\hbox{Shapes}(K)$
is finite, and $K$ is homeomorphic to $\hbox{\bb R}^n$, $n\ge 2$. If every
triangle in $K$ is isometric to its comparison triangle in $\hbox{\bb
R}^2$, then $K$ is isometric to $\hbox{\bb R}^n$.
\par
\medskip

\noindent {\bf Proof.}
It clearly suffices to show this for a neighbourhood of every vertex in
$K$.

Embed $n+1$ points into $K=K^n$ in general position, so that the chosen
vertex lies in the convex hull. The hull, $\Delta^n$, is homeomorphic to a
standard $n$-simplex. We show how to construct an isometry from this region
of $K$ to a region of $\hbox{\bb R}^n$. We construct this isometry on
successive skeleta of $\Delta^n$.

Consider a geodesic triangle forming part of the 1-skeleton of $\Delta^n$,
and consider a comparison triangle in $\hbox{\bb R}^n$. We automatically
have an isometry between these triangles.

Next we $\lq$fill-in' the triangles, that is, extend the isometry across
the interior. To do this, consider the one parameter family of geodesics
from a fixed vertex $v$ of the triangle in $K$ to each point of the
opposite side. The corresponding lines in our Euclidean triangle lead to a
one-to-one extension map in the obvious way. We claim that this is actually
an isometry.

Given any two points $x$, $y$, in the interior of the triangle in $K$, let
$\bar{x}$ and $\bar{y}$ denote the points where the extension of the
geodesics $vx$ respectively $vy$ meet the side of the original triangle
opposite to $v$. Our comparison triangle assumption applied to the triangle
$v\bar{x}\bar{y}$ means that the length of the side $x\bar{y}$ is the same
as the corresponding length in the Euclidean comparison triangle. Now apply
the assumption to the triangle $vx\bar{y}$ to show that the length of $xy$
is the same as the Euclidean comparison distance. Thus distances agree
under this mapping, as claimed.

If $n=2$, we are done. Otherwise, the original triangle in $K$ is a face of
some tetrahedron in $\Delta^n$. This tetrahedron has one vertex which is
not a vertex of the triangle. Consider another face of the tetrahedron. We
can find a comparison triangle in Euclidean space which intersects our
original comparison triangle in a common side. Again we can $\lq$fill-in'
and extend our mapping to the union of the two triangles.

By pivoting the second Euclidean triangle about the common edge, we can
clearly arrange for the distance between the two vertices not on the common
edge to be the same as that for the corresponding vertices in $\Delta^n$.
We then have that all faces of the tetrahedron defined by the points
introduced into $\hbox{\bb R}^n$ are comparison triangles for the
corresponding faces in $\Delta^n$. We can therefore extend our isometry to
an isometry of each face.

We claim that this is actually a global isometry of the union of the faces.
Given two points in different faces of $\Delta^n$, consider the minimal
geodesic in $\Delta^n$ linking them, and in particular consider the
point(s) at which the geodesic switches faces. Joining the corresponding
three (or more) points in the Euclidean picture, we obtain a curve of the
same length in $\hbox{\bb R}^n$. Suppose this is not a minimal geodesic for
our Euclidean tetrahedron: then the pre-image in $\Delta^n$ of the geodesic
which {\it is} minimal is a curve in $\Delta^n$ joining the given points of
strictly shorter length than the minimal geodesic. As this is impossible,
we deduce that corresponding minimal geodesics must have the same length,
and therefore we have an isometry of $\partial \Delta^n$ with the
corresponding complex in $\hbox{\bb R}^n$.

We now $\lq$fill-in' the isometry across the interior of the tetrahedron.
Repeating this process dimension by dimension gives the desired isometry
between $\Delta^n$ and some region of $\hbox{\bb R}^n.$ \hfill{\square}
\medskip

\noindent {\bf Proof of Theorem B.} The complex $K$ is a $\hbox{CAT}(0)$
space by Theorem A, and by Lemma 4.2 we have that $\hbox{Inv}_p(K)$ is a
geodesic space. Therefore [BFW; Proposition 6.2] applies, with the
conclusion that every geodesic triangle in $K$ is flat, that is, isometric
to its comparison triangle in $\hbox{\bb R}^2.$ By Lemma 4.3, $K$ must be
isometric to $\hbox{\bb R}^n.$ \hfill{\square}
\bigskip\bigskip\bigskip

\centerline{\bf REFERENCES}
\bigskip

\item{[BBI]} D. Burago, Y. Burago, S. Ivanov, {\it A course in metric
    geometry}, Graduate Studies in Mathematics vol. 33, American
    Mathematical Society, (2001).
\smallskip

\item{[BFW]} S. M. Buckley, K. Falk, D. J. Wraith, {\it Ptolemaic spaces
    and $\hbox{CAT}(0)$}, Glasgow J. Math. {\bf 51} (2009), 301--314.
\smallskip

\item{[BH]} M. R. Bridson, A. Haefliger, {\it Metric spaces of non-positive
    curvature}, Springer Verlag, (1999).
\smallskip

\item{[BHX]} S. M. Buckley, D. Herron, X. Xie, {\it Metric space
    inversions, quasihyperbolic distance, and uniform spaces}, Indiana J.
    Math. {\bf 57} (2008), 837--890.
\smallskip

\item{[C]} I. Chavel, {\it Riemannian Geometry: A Modern Introduction},
    Cambridge Univ. Press (1993).
\smallskip

\item{[FLS]} T. Foertsch, A. Lytchak, V. Schroeder, {\it Nonpositive
    curvature and the Ptolemy inequality}  Int. Math. Res. Not. IMRN
    (2007), article ID rnm100, 15 pages.
\smallskip

\item{[K]} D. C. Kay, {\it Ptolemaic metric spaces and the characterization
    of geodesics by vanishing metric curvature}, PhD Thesis, Michigan State
    University (1963).

\bigskip\bigskip

{\sl%
\noindent S.M.~Buckley and D.J.~Wraith:\hfill\break%
Department of Mathematics\hfill\break%
National University of Ireland Maynooth, Maynooth, Co. Kildare, Ireland.\hfill\break%
e-mail: {\tt stephen.buckley@nuim.ie and david.wraith@nuim.ie}
\medskip

\noindent J.~McDougall:\hfill\break%
Department of Mathematics and Computer Science\hfill\break%
Colorado College, Colorado Springs, Colorado 80903, USA.\hfill\break%
e-mail: {\tt JMcDougall@ColoradoCollege.edu.}
\medskip%
}

\end